\documentclass[amstex,16pt, amssymb]{article}

\usepackage{mathtext}
\usepackage[cp1251]{inputenc}
\usepackage[T2A]{fontenc}
\usepackage[dvips]{graphicx}
\usepackage{amsmath}
\usepackage{amssymb}
\usepackage{amsxtra}
\usepackage{latexsym}
\usepackage{ifthen}
\usepackage{amsthm}

\numberwithin{equation}{section}

\begin{document}

\title{To the theory of semi-linear Beltrami equations}

\author{V. Gutlyanski\u\i{}, O. Nesmelova, V. Ryazanov, E. Yakubov}

\date{}

\maketitle

\begin{abstract}
The present paper is devoted to the study of semi-linear Beltrami
equations which are closely relevant to the corresponding
semi-linear Poisson type equations of mathematical physics on the
plane in anisotropic and inhomogeneous media.

In its first part, applying completely continuous ope\-ra\-tors by
Ahlfors-Bers and Leray--Schauder approach, we prove existence of
regular solutions of the semi-linear Beltrami equations with no
boundary conditions. Moreover, here we derive their representation
through solutions of the Vekua type equations and generalized
analytic functions with sources.

As consequences, it is given a series of applications of these
results to semi-linear Poisson type equations and to the
corresponding equations of ma\-the\-ma\-ti\-cal physics describing
such phenomena as diffusion with phy\-si\-cal and chemical
absorption, plasma states and stationary burning in anisotropic and
inhomogeneous media.

The second part of the paper contains existence, representation and
regularity results for nonclassical solutions to the Hilbert
(Dirichlet) boundary value problem for semi-linear Beltrami
equations and to the Poincare (Neumann) boundary value problem for
semi-linear Poisson type equations with arbitrary boundary data that
are measurable with respect to logarithmic capacity.
\end{abstract}




\normalsize \baselineskip=18.5pt

\vskip 1cm


\section{Introduction}

Recall that the {\bf Beltrami equation} is the equation of the form
\begin{equation}\label{1}
f_{\bar{z}}=\mu(z) f_z
\end{equation} a.e. in $D$, where $\mu:
D\to\mathbb C$ is a mea\-su\-rab\-le function with $|\mu(z)|<1$
a.e.,  $f_{\bar z}=(f_x+if_y)/2$, $f_{z}=(f_x-if_y)/2$, $z=x+iy$,
$f_x$ and $f_y$ are partial derivatives of the function $f$ in $x$
and $y$, respectively. Note that continuous functions with
generalized derivative $f_{\bar z}=0$ are analytic functions, see
e.g. Lemma 1 in \cite{ABe}.


Equation~\eqref{1} is said to be {\bf nondegenerate} if
$||\mu||_{\infty}<1$, i.e., if $K_{\mu}\in L^{\infty}$,
\begin{equation}\label{K}
K_{\mu}(z)\ :=\ \frac{1\ +\ |\mu(z)|}{1\ -\ |\mu(z)|}\ .
\end{equation}
Homeomorphic solutions $f$ of nondegenerate (\ref{1}) in
$W^{1,2}_{\rm loc}$ are called {\bf quasiconformal mappings} or
sometimes {\bf $\mu -$conformal mappings}. Its continuous solutions
in $W^{1,2}_{\rm loc}$ are called {\bf $\mu -$conformal functions}.
On the cor\-res\-pon\-ding existence theorems for nondegenerate
(\ref{1}), see e.g. \cite{Alf}, \cite{BGMR} and \cite{LV}.

The {\bf inhomogeneous Beltrami equations}
\begin{equation}\label{s}
\omega_{\bar{z}}\ =\ \mu(z)\cdot \omega_z\ +\ \sigma(z)
\end{equation}
have been introduced and investigated in paper \cite{ABe}, see also
monograph \cite{Alf}. Boundary value problems for (\ref{s}) have
been researched in our preprint \cite{GNRY*}.

The present paper is devoted to the {\bf semi-linear Beltrami
equations} of the form
\begin{equation} \label{SB}
\omega_{\bar{z}}\ =\ \mu(z)\cdot \omega_z\ +\ \sigma(z)\cdot
q(\omega(z))\ ,
\end{equation}
where $\sigma :D\to\mathbb C$ belongs to class $L_p(D)$, $p>2,$
$q:\mathbb C\to\mathbb C$ is continuous and
\begin{equation} \label{q}
\lim\limits_{w\to \infty}\ \frac{q(w)}{w}\ =\ 0\ .
\end{equation}

For the history of semi-linear equations, see e.g. the monograph
\cite{MV} and our papers \cite{GNR-}--\cite{GNRY+} and \cite{R8}. We
prove here the existence of regular solutions of equations
(\ref{SB}) and derive their representation through solutions of the
Vekua type equations and generalized analytic functions with
sources.

In this connection, recall that the known monograph \cite{Ve} was
devoted to the {\bf generalized analytic functions}, i.e.,
continuous complex valued functions $H(z)$ of one complex variable
$z=x+iy$ of class $W^{1,1}_{\rm loc}$ satisfying the equations
\begin{equation}\label{eqG}
\partial_{\bar z}H\ +\ aH\ +\ b {\overline H}\ =\ s\ ,\ \ \ \ \ \partial_{\bar
z}\ :=\ \frac{1}{2}\left(\ \frac{\partial}{\partial x}\ +\
i\cdot\frac{\partial}{\partial y}\ \right)\ ,
\end{equation}
where it was assumed that the complex valued functions $a,b$ and $s$
belong to class $L^{p}(D)$ with some $p>2$ in the corresponding
domain $D\subseteq \mathbb C$.

Our papers \cite{GNRY} and \cite{R7} were devoted to {\bf
generalized analytic functions H with sources s}, when $a\equiv
0\equiv b$,
\begin{equation}\label{eqS}
\partial_{\bar z}H(z)\ =\ s(z)\ ,
\end{equation}
and the papers \cite{GNRY+} and \cite{R8} to the semi-linear Vekua
type equations of the form
\begin{equation} \label{SV}
\partial_{\bar z}H(z)\ =\ g(z)\cdot q(H(z))\
.
\end{equation}

Here we give a connection of the semi-linear Beltrami equations
(\ref{SB}) and {\bf semi-linear Poisson type equations} of the form
\begin{equation}\label{SD} {\rm div}\,
A(z)\,\nabla\,U(z)\ =\ G(z)\cdot Q(U(z)) \ ,
\end{equation} where $A(z)$ is a matrix function
associated with $\mu$ in (\ref{SB}), $G:D\to\mathbb R$ is a
multiplier in class $L_{p^{\prime}}(D)$, $p^{\prime}>1,$ and
$Q:\mathbb R\to\mathbb R$ is continuous such that
\begin{equation} \label{Q}
\lim\limits_{t\to \infty}\ \frac{Q(t)}{t}\ =\ 0\ .
\end{equation}

As a consequence, we obtain also the existence of weak regular
solutions of equations (\ref{SD}). Moreover, we derive the
representation of the given solutions of equations (\ref{SD})
through regular solutions $h$ of semi-linear Poisson equations
\begin{equation}\label{SP} \triangle\,h(z)\ =\ g(z)\cdot Q(h(z))
\end{equation}
which are the so--called weak generalized harmonic functions with
sources.


Recall that, correspondingly to \cite{GNRY}, \cite{R7}, a continuous
function $h:D\to\mathbb R$ of class $W^{2,p}_{\rm loc}$ is called a
{\bf generalized harmonic function with a source s}:$\,D\to\mathbb
R$ in $L_p(D)$, $p>2$, if $h$ a.e. satisfies the Poisson equation
\begin{equation}\label{eqSSS}
\triangle h(z)\ =\ s(z)\ .
\end{equation} Note that by the Sobolev embedding theorem, see
Theorem I.10.2 in \cite{So}, such functions $h$ belong to the class
$C^1$.

Here we say that a function $h:D\to\mathbb R$ is a {\bf weak
generalized harmonic function with a source s} $:D\to\mathbb C$ in
$L_p(D)$, $p>2$, if $h$ is a real part of a generalized analytic
function with the source $s$. We show that every generalized
harmonic function with a real source $s:D\to\mathbb R$ is a weak
generalized harmonic function with the source $s$, see Remark 5
below. However, it is clear that the inverse conclusion has no sense
for complex sources $s$ at all.


The equations (\ref{SD}) describe many physical phenomena as
diffusion with  physical and chemical absorption, plasma states,
stationary burning etc. in anisotropic and inhomogeneous media.
Hence we have been able to give the corresponding applications of
the obtained results to the mathematical physics.

\bigskip

\section{Factorization of semi-linear Beltrami equations}

Later on, we always assume that, for the corresponding $p>2$,
\begin{equation}\label{b}
 k\,C_p\ <\ 1\ ,\ \ \ \ k\ :=\ \|\mu\|_{\infty}\ <\ 1\ ,
\end{equation}
where $C_p$ is the norm of the known operator $T:L_p(\mathbb C)\to
L_p(\mathbb C)$ defined through the Cauchy principal limit of the
singular integral
\begin{equation}\label{i}
 (Tg)(\zeta)\ :=\
 \lim\limits_{\varepsilon\to 0}\left\{-\frac{1}{\pi}\int\limits_{|z-\zeta|>\varepsilon}\frac{g(z)}{(z-\zeta)^2}\
 dxdy\right\}\ ,\ \ \ z=x+iy\ .
\end{equation}

\medskip

The next lemmas gives the factorization of solutions of semi-linear
Beltrami equations (\ref{SB}) through solutions of  semi-linear
Vekua type equations (\ref{SV}).

\medskip

{\bf Lemma 1.} {\it Let $D$ be a bounded domain in $\mathbb C$, $\mu
: D\to\mathbb C$ be measurable with $\|\mu\|_{\infty}<1$,
$\sigma:D\to\mathbb C$ be in class $L_p(D)$, $p>2$. Suppose that
$q:\mathbb C\to\mathbb C$ is continuous and $f^{\mu}:\mathbb
C\to\mathbb C$ is a $\mu-$conformal mapping from Theorem B in
\cite{GNRY*} with an extension of $\mu$ onto $\mathbb C$ keeping
compact support and condition (\ref{b}).

Then each continuous solution $\omega$ of equation (\ref{SB}) in $D$
of class $W^{1,p}(D)$ has the representation as a composition
$H\circ f^{\mu}|_D$, where $H$ is a continuous solution of
(\ref{SV}) in class  $W^{1,p_*}_{\rm loc}(D_*)$, where
$D_*:=f^{\mu}(D)$, $p_*:=p^2/2(p-1)\in(2,p)$, with the multiplier
$g$ in (\ref{SV}) of class $L_{p_*}(D_*)$ defined by formula
\begin{equation}\label{r} g\ :=\
\left(\frac{f^{\mu}_z}{J^{\mu}}\cdot\sigma\right)\circ
\left(f^{\mu}\right)^{-1}\ ,
\end{equation} where $J^{\mu}$ is the Jacobian of $f^{\mu}$.

Vice versa, if $H$ is a continuous solution in class $W^{1,p_*}_{\rm
loc}(D_*)$ of (\ref{SV}) with multiplier $g\in L_{p_*}(D_*)$,
$p_*>2$, given by (\ref{r}), then $\omega :=H\circ f^{\mu}$ is a
solution of (\ref{SB}) in class $C^{\alpha}_{\rm loc}\cap
W^{1,p^*}_{\rm loc}(D)$, $\alpha = 1-2/p^*$, where
$p^*:=p_*^2/2(p_*-1)\in(2,p_*)$.}

\begin{proof}
Indeed, if $\omega$ is a continuous solution of (\ref{SB}) in $D$ of
class $W^{1,p}(D)$, then $\omega$ is a solution of (\ref{s}) in $D$
with the source $\Sigma :=\sigma\cdot q\circ\omega$ in the same
class. Then by Lemma 1 and Remark 1 in \cite{GNRY*} $\omega = H\circ
f^{\mu}$, where $H$ is a generalized analytic function with the
source $G$ of class $L_{p_*}(D_*)$ after replacement of $\sigma$ by
$\Sigma$ in (\ref{r}). Note that $H\in W^{1,p_*}_{\rm loc}(D_*)$,
see e.g. Theorems 1.16 and 1.37 in \cite{Ve}. The proof of the vice
versa conclusion of Lemma 2 is similar and it is again based on its
reduction to Lemma 1 in \cite{GNRY*}.
\end{proof}


\section{On solutions of semi-linear Beltrami equations}

First of all, recall that a {\bf com\-ple\-te\-ly continuous}
mapping from a metric space $M_1$ into a metric space $M_2$ is
defined as a continuous mapping on $M_1$ which takes bounded subsets
of $M_1$ into relatively compact subsets of $M_2$, i.e., with
compact closures in space $M_2$. When a continuous mapping takes
$M_1$ into a relatively compact subset of $M_1$, it is nowadays said
to be {\bf compact} on $M_1$.

\medskip

Note that the notion of completely continuous (compact) operators is
due essentially to Hilbert in a special space that, in reflexive
spaces, is equivalent to Definition VI.5.1 for the Banach spaces in
\cite{DS} which is due to F. Riesz, see also further comments of
Section VI.12 in \cite{DS}.

\medskip

Recall some further definitions and the fundamental result of the
celebrated paper \cite{LS}. Leray and Schauder extend as follows the
Brouwer degree to compact perturbations of the identity $I$ in a
Banach space $B$, i.e., a complete normed linear space. Namely,
given an open bounded set $\Omega\subset B$, a compact mapping $F:
B\to B$ and $z \notin \Phi(\partial \Omega)$, $\Phi :=I-F$, the {\bf
(Leray–Schauder) topological degree} $\deg\, [\Phi,\Omega, z]$ of
$\Phi$ in $\Omega$ over $z$ is constructed from the Brouwer degree
by approximating the mapping $F$ over $\Omega$ by mappings
$F_{\varepsilon}$ with range in a finite-dimensional subspace
$B_{\varepsilon}$ (containing $z$) of $B$. It is showing that the
Brouwer degrees $\deg\, [\Phi_{\varepsilon} ,\Omega_{\varepsilon},
z]$ of $\Phi_{\varepsilon}:=I_{\varepsilon} - F_{\varepsilon}$,
$I_{\varepsilon}:=I|_{B_{\varepsilon}}$, in
$\Omega_{\varepsilon}:=\Omega\cap B_{\varepsilon}$ over $z$
stabilize for sufficiently small positive $\varepsilon$ to a common
value defining $\deg\, [\Phi,\Omega, z]$ of $\Phi$ in $\Omega$ over
$z$.

\medskip

This topological degree “algebraically counts” the number of fixed
points of $F(\cdot)-z$ in $\Omega$ and conserves the basic
properties of the Brouwer degree as ad\-di\-ti\-vi\-ty and homotopy
invariance. Now, let $a$ be an isolated fixed point of $F$. Then the
{\bf local (Leray–Schauder) index} of $a$ is defined by ${\rm ind}\,
[\Phi, a] := \deg [\Phi,B(a, r), 0]$ for small enough $r > 0$. ${\rm
ind}\, [\Phi, 0]$ is called by {\bf index} of $F$. In particular, if
$F\equiv 0$, correspondingly, $\Phi\equiv I$, then the index of $F$
is equal to $1$.

\medskip

Let us formulate the main result in \cite{LS}, Theorem 1, see also
the survey \cite{Ma}.

\medskip

{\bf Proposition 1.} {\it Let $B$ be a Banach space, and let
$F(\cdot,\tau):B\to B$ be a family of operators with $\tau\in[0,1]$.
Suppose that the following hypotheses hold:


{\rm {\bf (H1)}} $F(\cdot,\tau)$ is completely continuous on $B$ for
each $\tau\in[0,1]$ and uniformly continuous with respect to the
parameter $\tau\in[0,1]$ on each bounded set in $B$;

{\rm {\bf (H2)}} the operator $F:=F(\cdot,0)$ has finite collection
of fixed points whose total index is not equal to zero;

{\rm {\bf (H3)}} the collection of all fixed points of the operators
$F(\cdot,\tau)$, $\tau\in[0,1]$,  is bounded in $B$.

Then the collection of all fixed points of the family of operators
$F(\cdot,\tau)$ contains a continuum along which $\tau$ takes all
values in  $[0,1]$.}

\medskip

We denote by $B_p$ the Banach space of functions $\omega:\mathbb
C\to\mathbb C$, which satisfy a global H\"older condition of order
$1 - 2/p$, which vanish at the origin, and whose generalized
derivatives $\omega_z$ and $\omega_{\bar z}$ exist and belong to
$L_p(\mathbb C)$. The norm in $B_p$ is defined by
\begin{equation}\label{n}
\|\omega\|_{B_p}\ :=\ \sup\limits_{\underset{z_1\ne
z_2}{z_1,z_2\in\mathbb
C,}}\frac{|\omega(z_1)-\omega(z_2)|}{|z_1-z_2|^{1-2/p}}\ +\
\|\omega_z\|_p\ +\ \|\omega_{\bar z}\|_p\ .
\end{equation}

{\bf Remark 1.} By Lemma 5 in \cite{ABe} the mapping
$\sigma\to\omega^{\mu,\sigma}$ from Theorem A in \cite{GNRY*} is a
bounded linear operator from $L_p(\mathbb C)$ to $B_p(\mathbb C)$
with a bound that depends only on $k$ and $p$ in (\ref{b}). In
particular, this is a bounded linear operator from $L_p(\mathbb C)$
to $C(\mathbb C)$. Namely, by (15) in \cite{ABe} we have that
$\omega^{\mu,\sigma}$ is H\"older continuous:
\begin{equation}\label{H}
|\omega^{\mu,\sigma}(z_1)\ -\ \omega^{\mu,\sigma}(z_2)|\ \leq\
c\cdot\|\sigma\|_p\cdot |z_1 - z_2|^{1-2/p}\ \ \ \ \ \forall\ z_1\,\
z_2\ \in \mathbb C\ ,
\end{equation}
where the constant $c$ may depend only on $k$ and $p$ in (\ref{b}).
Moreover, $\omega^{\mu,\sigma}(z)$ is locally bounded because
$\omega^{\mu,\sigma}(0)=0$. Thus, the linear operator
$\sigma\to\omega^{\mu,\sigma}|_S$ is completely continuous for each
compact set $S$ in $\mathbb C$ by Arzela-Ascoli theorem, see e.g.
Theorem IV.6.7 in \cite{DS}.

\medskip

{\bf Theorem 1.} {\it Let $\mu :\mathbb C\to\mathbb C$ and $\sigma
:\mathbb C\to\mathbb C$ have compact supports, $\mu\in
L_{\infty}(\mathbb C)$ with $k:=\|\mu\|_{\infty}<1$ and $\sigma \in
L_p(\mathbb C)$ for some $p>2$ satisfying (\ref{b}). Suppose that
$q:\mathbb C\to\mathbb C$ is a continuous function satisfying
condition (\ref{q}). Then the semi-linear Beltrami equation
(\ref{SB}) has a solution $\omega$ of class $B_p(\mathbb C)$.}

\medskip

{\bf Remark 2.} By Lemma 1 $\omega$ has the representation as a
composition $H\circ f^{\mu}$, where $f^{\mu}:\mathbb C\to\mathbb C$
is a $\mu-$conformal mapping from Theorem B in \cite{GNRY*} and $H$
is a continuous solution of (\ref{SV}) in class $W^{1,p_*}_{\rm
loc}(\mathbb C)$, $p_*:=p^2/2(p-1)\in(2,p)$, with the multiplier $g$
in (\ref{SV}) of class $L_{p_*}(\mathbb C)$ defined by formula
(\ref{r}). Note also that $H$ is a generalized analytic function
with a source of the same class.

\begin{proof} If $\| \sigma\|_p=0$ or $\| q\|_C= 0$, then Theorem A
in \cite{GNRY*} above gives the desired solution $\omega
:=\omega^{\mu , 0}$ of equation (\ref{SB}). Thus, we may assume that
$\| \sigma\|_p\neq 0$ and $\| q\|_C\neq 0$. Set
$q_*(t)=\max\limits_{|w|\le t}|q(w)|$, $t\in\mathbb
R^+:=[0,\infty)$. Then the function $q_*:\mathbb R^+\to\mathbb R^+$
is continuous and nondecreasing and, moreover, by (\ref{q})
\begin{equation} \label{*}
\lim\limits_{t\to \infty}\ \frac{q_*(t)}{t}\ =\ 0\ .
\end{equation}

Let us show that the family of operators $F(g;\tau):
L^{\sigma}_{p}(\mathbb C)\to L^{\sigma}_{p}(\mathbb C)$,
\begin{equation} \label{eqFORMULA}
F(g;\tau)\ :=\ \tau \sigma\cdot q(\omega^{\mu , g}) \ \ \ \ \ \ \
\forall\ \tau\in[0,1]\ ,
\end{equation}
where $L^{\sigma}_{p}(\mathbb C)$ consists of functions $g\in
L_{p}(\mathbb C)$ with supports in the support $S$ of $\sigma$,
satisfies hypotheses H1-H3 of Theorem 1 in \cite{LS}, see
Proposition 1 above. Indeed:

H1). First of all, the function $F(g;\tau)\in L^{\sigma}_{p}(\mathbb
C)$ for all $\tau\in[0,1]$ and  $g\in L^{\sigma}_{p}(\mathbb C)$
because the function $q(\omega^{\mu ,g})$ is continuous and,
furthermore, the operators $F(\cdot ;\tau)$ are completely
continuous for each $\tau\in[0,1]$ and even uniformly continuous
with respect to parameter $\tau\in[0,1]$ by Theorem A in
\cite{GNRY*} and Remark 1.

H2). The index of the operator $F(g;0)$ is obviously equal to $1$.

H3). Let us assume that the collection of all solutions of the
equations $g=F(g;\tau)$, $\tau\in[0,1]$, is not bounded in
$L^{\sigma}_{p}(\mathbb C)$, i.e., there is a sequence of functions
$g_n\in L^{\sigma}_{p}(\mathbb C)$ with $\|g_n\|_p\to\infty$ as
$n\to\infty$ such that $g_n=F(g_n;\tau_n)$ for some
$\tau_n\in[0,1]$, $n=1,2,\ldots$.

However, then by Remark 1 we have that
$$
\| g_n\|_p\ \le\ \| \sigma\|_p\ q_*\left(\| \omega^{\mu ,
g_n}|_S\|_C\right)\ \le\ \| \sigma\|_p\ q_*\left(\, M\, \|
g_n\|_p\right)
$$
for some constant $M>0$ and, consequently,
\begin{equation} \label{eqEST3W}
\frac{q_*(\, M\, \| g_n\|_p)}{M\, \| g_n\|_p}\ \ge\ \frac{1}{M\, \|
\sigma\|_p}\ >\ 0\ .
\end{equation}
The latter is impossible by condition (\ref{*}). The obtained
contradiction disproves the above assumption.

Thus, by Theorem 1 in \cite{LS}, see Proposition 1 above, there is a
function $g\in L^{\sigma}_p(\mathbb C)$ with $F(g;1)=g$, and then by
Theorem A in \cite{GNRY*} the function $\omega:=\omega^{\mu , g}$
gives the desired solution of (\ref{SB}).
\end{proof}

Let us also give the following lemma on semi-linear Beltrami
equations that is of independent interest, in particular, for
significant applications to semi-linear analogs of the Poisson
equations in anisotropic and inhomogeneous media.

\medskip

{\bf Lemma 2.} {\it Let $D$ be a bounded domain in $\mathbb C$, $\mu
:D\to\mathbb C$ in class $L_{\infty}(D)$, $k:=\|\mu\|_{\infty}<1$,
$G : D\to\mathbb C$ in $L_p(D)$ for some $p^{\prime}>1$ and ${\cal
L}:L_{p^{\prime}}(D)\to L_p(D)$ be a linear bounded operator for
some $p>2$ satisfying (\ref{b}).

Suppose that $q:\mathbb C\to\mathbb C$ is a continuous function
satisfying condition (\ref{q}). Then the semi-linear Beltrami
equation of the form
\begin{equation} \label{SBL}
\omega_{\bar{z}}\ =\ \mu(z)\cdot \omega_z\ +\ {\cal L}[G
q(\omega)](z)\ ,\ \ \ \ z\in D\ ,
\end{equation}
has a solution $\omega$ of class $C^{\alpha}(D)\cap W^{1,p}(D)$ with
$\alpha =1-2/p$.}

\begin{proof}
Indeed, arguing perfectly similar to the proof of Theorem 1 for
\begin{equation} \label{eqFORMULAL}
F(g;\tau)\ :=\ {\cal L}\left[\tau G q(\omega^{\mu , g})\right]\ :\
L_{p}(D)\to L_{p}(D) \ , \ \ \ \ \ \ \ \tau\in[0,1]
\end{equation} with $\mu$, $G$, $g$ and ${\cal L}$ extended by zero outside of $D$,
we see that the family of the operators $F(g;\tau)$, $\tau\in[0,1]$,
satisfies all the hypotheses of Theorem 1 in \cite{LS}, see
Proposition 1 above. Thus, there is $g\in L_p(\mathbb C)$ with
$F(g;1)=g$, and then by Theorem A in \cite{GNRY*} the function
$\omega:=\omega^{\mu , g}|_D$ gives the desired solution of
(\ref{SBL}).
\end{proof}

{\bf Remark 3.} Moreover, arguing similarly to the proof of Lemma 1
one can show that $\omega =H\circ f^{\mu}|_D$, where
$f^{\mu}:\mathbb C\to\mathbb C$ is a $\mu-$conformal mapping from
Theorem B in \cite{GNRY*} with $\mu$ extended onto $\mathbb C$ by
zero outside of $D$ and $H:D_*\to\mathbb C$ is a generalized
analytic function in the domain $D_*:=f^{\mu}(D)$ with the source
\begin{equation}\label{rL}
s\ :=\ \left\{\frac{f^{\mu}_z}{J^{\mu}}\cdot{\cal L}\left[G
q(\omega)\right]\right\}\circ \left(f^{\mu}\right)^{-1} \in\
L_{p_*}(D_*)\ ,
\end{equation} where $J^{\mu}$ is the Jacobian of $f^{\mu}$ and $p_*:=p^2/2(p-1)\in(2,p)$.

\section{Toward semi-linear Poisson type equations}

Further $\mathbb S^{2\times 2}$ denotes the collection of all
$2\times 2$ matrices with real elements
\begin{equation}\label{matrix}
A\  =\ \left[\begin{array}{ccc} a_{11}  & a_{12} \\
a_{21} & a_{22} \end{array}\right]\end{equation} which are
symmetric, i.e., $a_{12}=a_{21}$, with ${\rm det}\,A=1$ and {\bf
ellipticity condition} ${\rm det}\,(I+A)>0$, where $I$ is the unit
$2\times 2$ matrix. The latter condition means in terms of elements
of $A$ that $(1+a_{11})(1+a_{22})>a_{12}a_{21}$.

Now, let us first consider in a domain $D$ of the complex plane
$\mathbb C$ the linear Poisson type equations
\begin{equation}\label{6.8} {\rm div}\,
[A(z)\,\nabla\,u(z)]\ =\ g(z)\ , \end{equation} where $A: D\to
\mathbb S^{2\times 2}$ is a measurable matrix function  whose
elements $a_{ij}(z)$, $i,j=1,2$ are bounded and $g:D\to\mathbb R$ is
a scalar function in $L^1_{\rm loc}$.

Note that (\ref{6.8}) are the main equation of hydromechanics
(mechanics of incompressible fluids) in anisotropic and
inhomogeneous media.

We say that a function $u:D\to\mathbb R$ is a {\bf generalized ${\bf
A-}$harmonic function with the source g}, see e.g. monograph
\cite{HKM}, if $u$ is a weak solution of (\ref{6.8}), i.e., if $u\in
W^{1,1}_{\rm loc}(D)$ and
\begin{equation}\label{weak}
\int\limits_D \langle A(z)\nabla u(z),\nabla\psi(z)\rangle\, d\,m(z)
+ \int\limits_D g(z)\,\psi(z)\ d\,m(z) =\ 0
\end{equation}
for all $\psi\in C_0^\infty(D)$, where $C_0^\infty(D)$ denotes the
collection of all infinitely differentiable functions $\psi:
D\to\mathbb R$ with compact support in $D$, $\langle a,b\rangle$
means the scalar product of vectors $a$ and $b$ in $\mathbb R^2$,
and $d\,m(z)$ corresponds to the Lebesgue measure  in the plane
$\mathbb C$.

Later on, we use the {\bf logarithmic (Newtonian) potential of
sources} $g\in L^1(\mathbb C)$ with compact supports given by the
formula:
\begin{equation} \label{eqIPOTENTIAL} {\cal N}^{g}(z)\  :=\
\frac{1}{2\pi}\int\limits_{\mathbb C} \ln|z-w|\, g(w)\ d\, m(w)\
.\end{equation}

By Lemma 3 in \cite{GNR1} and in \cite{GNR2} we have its following
basic properties.


{\bf Remark 4.}  Let $g:\mathbb C\to\mathbb R$ have compact support.
If $g\in L^1(\mathbb C)$, then ${\cal N}^g\in L^{r}_{\rm
loc}(\mathbb C)$ for all $r\in [1,\infty)$, ${\cal N}^g\in
W^{1,p}_{\rm loc}(\mathbb C)$ for all $p\in [1,2)$, moreover, there
exist generalized derivatives by Sobolev $\frac{\partial^2
N^{g}}{\partial z\partial\overline z}$ and $\frac{\partial^2
N^{g}}{\partial\overline z\partial z}$ satisfying the equalities
\begin{equation}\label{eqLAP} 4\cdot\frac{\partial^2 N^{g}}{\partial
z\partial\overline z}\ =\ \triangle N^g\ =\ 4\cdot\frac{\partial^2
N^{g}}{\partial\overline z\partial z}\ =\ g\ \ \ \mbox{a.e.}
\end{equation}
Furthermore, if $g\in L^{p^{\prime}}(\mathbb C)$ for some
$p^{\prime}>1$, then ${\cal N}^g\in W^{2,p^{\prime}}_{\rm
loc}(\mathbb C)$, moreover, ${\cal N}^g\in W^{1,p}_{\rm loc}(\mathbb
C)$ for some $p>2$ and, consequently, ${\cal N}^g\in C^{\alpha}_{\rm
loc}(\mathbb C)$ with $\alpha =1-2/p$. Finally, if $g\in
L^{p^{\prime}}(\mathbb C)$ for some $p^{\prime}>2$, then  ${\cal
N}^g\in C^{1,\alpha}_{\rm loc}(\mathbb C)$ with $\alpha =1-2/p$.


Next, we say that a function $v:D\to\mathbb R$ is {\bf ${\bf
A-}$conjugate of a ge\-ne\-ra\-li\-zed ${\bf A-}$harmonic function u
with a source g}$:D\to\mathbb R$ if $v\in W^{1,1}_{\rm loc}(D)$ and
\begin{equation}\label{eqAconjugate}
\nabla v(z)\ =\ {\mathbb H} \left[\, A(z)\nabla u(z)\ -\ \nabla
{\cal N}^g(z)\, \right]\ \ \ \mbox{a.e.}\ ,\ \ \ \ \ \  {\mathbb H}\ :=\ \left[\begin{array}{rr} 0 &-1 \\
        1 &0\end{array}\right]\ .
\end{equation}


{\bf Lemma 3.} {\it Let $D$ be a bounded domain in $\mathbb C$,
$g:D\to\mathbb R$ be in $L^1(D)$ and let $u$ be a weak solution of
equation (\ref{6.8}) with a matrix function $A:D\to\mathbb
S^{2\times 2}$ whose elements $a_{ij}(z)$, $i,j=1,2$ are bounded and
measurable.

If $v$ is $A-$conjugate of $u$, then $\omega :=u+iv$ satisfies the
nondegenerate in\-ho\-mo\-geneous Beltrami equation (\ref{s}) with
\begin{equation}\label{6.10}
\mu(z)\ :=\ \mu_A(z)\ =\frac{1}{\mathrm{det}\, \left[\,
I+A(z)\,\right]}\ \left[\,a_{22}(z)-a_{11}(z)\ -\
2ia_{21}(z)\,\right]\ ,\end{equation}
\begin{equation}\label{sigma}
\sigma(z)\ :=\ {\cal N}_{\bar z}^g(z)\ +\ \mu(z)\, \overline{{\cal
N}_{\bar z}^g(z)}\ .
\end{equation}

Conversely, if $\omega\in W^{1,1}_{\rm loc}(D)$ is a solution of the
nondegenerate in\-ho\-mo\-ge\-ne\-ous Beltrami equation (\ref{s})
with $\sigma$ given by (\ref{sigma}), then $u:={\rm Re}\,\omega$ is
a weak solution of equation (\ref{6.8}) with the matrix valued
function $A:D\to\mathbb S^{2\times 2}$,
\begin{equation}\label{6.9}
A(z)\ :=\ \left[\begin{array}{ccc} {|1-\mu(z)|^2\over 1-|\mu(z)|^2}  & {-2{\rm Im}\,\mu(z)\over 1-|\mu(z)|^2} \\
                            {-2{\rm Im}\,\mu(z)\over 1-|\mu(z)|^2}          & {|1+\mu(z)|^2\over 1-|\mu(z)|^2}
                            \end{array}\right]\ ,
\end{equation} whose elements are bounded and
measurable.}


{\bf Remark 5.} Hence, in the case $A\equiv I$ and $g\in
L^{p^{\prime}}(D)$, $p^{\prime}>2$, we conclude that every
generalized harmonic function $u$ with the source $g$ is a weak
ge\-ne\-ra\-li\-zed harmonic function with the same source, see e.g.
Theorem 1.16 in \cite{Ve}. The inverse conclusion is, generally
speaking, not true and has no sense at all.

\begin{proof}
Indeed, let $u$ be a weak solution of equation (\ref{6.8}) with
$g:D\to\mathbb R$ in $L^1(D)$ and a matrix function $A:D\to\mathbb
S^{2\times 2}$ whose elements are bounded and measurable. Then by
(\ref{eqLAP}) we have that $u$ is a weak solution of the equation
\begin{equation}\label{As} {\rm div}\,
[A(z)\,\nabla\,u(z)]\ =\ {\rm div}\, \nabla {\cal N}^g(z)\ .
\end{equation}
If $v$ is $A-$conjugate of $u$, then by Theorem 16.1.6 in \cite{AIM}
the function $\omega :=u+iv$ satisfies the nondegenerate
in\-ho\-mo\-geneous Beltrami equation (\ref{s}) with $\mu$ and
$\sigma$ given by (\ref{6.10}) and (\ref{sigma}).

Conversely, if $\omega\in W^{1,1}_{\rm loc}(D)$ is a solution of the
nondegenerate in\-ho\-mo\-ge\-ne\-ous Beltrami equation (\ref{s})
with $\sigma$ given by (\ref{sigma}), then, again by Theorem 16.1.6
in \cite{AIM}, the functions $u:={\rm Re}\,\omega$ and $v:={\rm
Im}\,\omega$ satisfy the relation (\ref{eqAconjugate}) with the
matrix function $A:D\to\mathbb S^{2\times 2}$ given by (\ref{6.9})
whose elements $a_{ij}(z)$ are measurable in $z\in D$ and bounded
because $|a_{ij}|\leq \| K_{\mu}\|_{\infty}$. Note that
(\ref{eqAconjugate}) is equivalent to the equation
\begin{equation}\label{eqDIV}
A(z)\nabla\,u(z)\ -\ \nabla{\cal N}^{g}(z)\ =\ -{\mathbb H}\nabla\,
v(z)
\end{equation} because $\mathbb H^2=-I$.
As known, the curl of any gradient field is zero in the sense of
distributions and, moreover, the Hodge operator ${\mathbb H}$
transforms curl-free fields into divergence-free fields, and vice
versa, see e.g. 16.1.3 in \cite{AIM}. Hence $u$ is a weak solution
of equation (\ref{As}) as well as of (\ref{6.8}) in view of
(\ref{eqLAP}).
\end{proof}

{\bf Theorem 2}. {\it Let $D$ be a bounded domain in $\mathbb C$, a
scalar function $G:D\to\mathbb R$ be in class $L^{p^{\prime}}(D)$
for some $p^{\prime}>1$, a continuous function $Q:\mathbb
R\to\mathbb R$ satisfy condition (\ref{Q}) and let $A:D\to\mathbb
S^{2\times 2}$ be a matrix function whose elements $a_{ij}(z)$,
$i,j=1,2$ are bounded and measurable.

Then the semi-linear Poisson type equation (\ref{SD}) has a { weak
solution} $u$ of class $C^{\alpha}(D)\cap W^{1,p}(D)$ with $\alpha
=1-2/p$ for some $p>2$.

Moreover, $u=h\circ f^{\mu}|_D$, where $f^{\mu}:\mathbb C\to\mathbb
C$ is a $\mu-$conformal mapping from Theorem B and $h$ is a weak
generalized harmonic function in the domain $D_*:=f^{\mu}(D)$ with
the source
\begin{equation}\label{LrL}
s\ :=\ \left\{\frac{f^{\mu}_z}{J^{\mu}}\cdot\sigma\right\}\circ
\left(f^{\mu}\right)^{-1} \in\ L_{p_*}(D_*)\ ,
\end{equation} where $J^{\mu}$ is the Jacobian of $f^{\mu}$,
$p_*:=p^2/2(p-1)\in(2,p)$, $\mu$ is defined by formula (\ref{6.10})
and $\sigma$ is calculated by formula (\ref{sigma}) with
$g=G\,Q(u)$.}

\medskip

Here we say that a function $u:D\to\mathbb R$ is a {\bf weak
solution} of (\ref{SD}), if $u\in W^{1,1}_{\rm loc}(D)$ and
\begin{equation}\label{SEMIweak}
\int\limits_D \langle A(z)\nabla u(z),\nabla\psi(z)\rangle\,
d\,m(z)\ + \int\limits_D G(z)\,Q(u(z))\,\psi(z)\, d\,m(z)\ =\ 0
\end{equation}
for all $\psi\in C_0^\infty(D)$, where $C_0^\infty(D)$ denotes the
collection of all infinitely differentiable functions $\psi:
D\to\mathbb R$ with compact support in $D$, $\langle a,b\rangle$
means the scalar product of vectors $a$ and $b$ in $\mathbb R^2$,
and $d\,m(z)$ corresponds to the Lebesgue measure  in the plane
$\mathbb C$.

\begin{proof}
With no loss of generality, we may assume here that
$p^{\prime}\in(1,2]$ and that $g\equiv 0$ outside of $D$, and then
${\cal N}^g\in W^{1,p}(D)$ for all $p\in(1,p^*)$, where
$p^*=2p^{\prime}/(2-p^{\prime})>2$, see Lemma 3 in \cite{GNR1}.
Hence later on we may also assume that $p>2$ satisfies condition
(\ref{b}) for $\mu$ in (\ref{6.10}). Moreover, again by Lemma 3 in
\cite{GNR1}, the correspondence $g\to{\cal N}^g_{\bar z}$ generates
a completely continuous linear operator $L$ acting from real valued
$L_{p^{\prime}}(D)$ to complex valued $L_{p}(D)$. Thus, the linear
operator ${\cal L}:=L+\mu \bar L$ with the multiplier $\mu\in
L^{\infty}(D)$ is bounded. Then by Lemma 2 the semi-linear Beltrami
equation (\ref{SBL}) with $q(\omega):=Q({\rm Re}\,\omega)$ has a
solution $\omega$ of class $C^{\alpha}(D)\cap W^{1,p}(D)$ with
$\alpha =1-2/p$. Moreover, by Lemma 3 the function $u:={\rm
Re}\,\omega$ is a weak solution of equation (\ref{6.8}) of the given
class. Finally, by Lemma 1 in \cite{GNRY*} we conclude that $u$ has
the representation as the composition $h\circ f^{\mu}|_D$, where
$f^{\mu}:\mathbb C\to\mathbb C$ is a $\mu-$conformal mapping from
Theorem B and $h$ is a weak generalized harmonic function in the
domain $D_*:=f^{\mu}(D)$ with the source (\ref{LrL}).
\end{proof}

\section{Semi-linear Poisson type equations in physics}

Theorem 2 on semi-linear Poisson type equations can be applied to
ma\-the\-ma\-ti\-cal models of physical and chemical absorption with
diffusion, plasma states, stationary burning etc. in anisotropic and
inhomogeneous media.

The first group of such applications is relevant to
reaction-diffusion problems. Problems of this type are discussed in
\cite{Diaz}, p. 4, and, in details, in \cite{Aris}. A nonlinear
system is obtained for the density $U$ and the temperature $T$ of
the reactant. Upon eliminating $T$ the system can be reduced to
equations of the form
\begin{equation}\label{RDP}
 \triangle  U\ =\ \sigma\cdot Q(U)
\end{equation}
with $\sigma >0$ and, for isothermal reactions, $Q(U) = U^{\lambda}$
where $\lambda>0$ that is called the order of the reaction. It turns
out that the density of the reactant $U$ may be zero in a subdomain
called a {\bf dead core}. A particularization of results in Chapter
1 of \cite{Diaz} shows that a dead core may exist just if and only
if $\beta\in(0,1)$, see also the corresponding examples in
\cite{GNR}.

In the case of anisotropic and inhomogeneous media, we come to the
semi-linear Poisson type equations (\ref{SD}). In this connection,
the following statement may be of independent interest.

\medskip

{\bf Corollary 1.} {\it Let $D$ be a bounded domain in $\mathbb C$,
a scalar function $\sigma :D\to\mathbb R$ be in class
$L^{p^{\prime}}(D)$ for some $p^{\prime}>1$ and let $A:D\to\mathbb
S^{2\times 2}$ be a matrix function whose elements $a_{ij}(z)$,
$i,j=1,2$ are bounded and measurable.

Then there is a weak solution $u:D\to\mathbb R$ of class
$C^{\alpha}_{\rm loc}\cap W^{1,p}_{\rm loc}$ with $\alpha =1-2/p$
for some $p>2$ to the semi-linear Poisson type equation
\begin{equation} \label{eqQUASILINEARW568C}
{\rm div}\, [A(z)\,\nabla\, u(z)]\ =\ \sigma(z)\cdot u^{\lambda}(z)\
,\ \ \ 0\ <\ \lambda\ <\ 1\  ,\ \ \ \ \ \ \ \mbox{a.e. in $D$\ .}
\end{equation}}


Note also that certain mathematical models of a thermal evolution of
a heated plasma lead to nonlinear equations of the type (\ref{RDP}).
Indeed, it is known that some of them have the form
$\triangle\psi(u)=f(u)$ with $\psi'(0)=\infty$ and $\psi'(u)>0$ if
$u\not=0$ as, for instance, $\psi(u)=|u|^{q-1}u$ under $0 < q < 1$,
see e.g. \cite{Diaz}. With the replacement of the function
$U=\psi(u)=|u|^q\cdot {\rm sign}\, u$, we have that $u = |U|^Q\cdot
{\rm sign}\, U$, $Q=1/q$, and, with the choice $f(u) =
|u|^{q^2}\cdot {\rm sign}\, u$, we come to the equation $\triangle U
= |U|^q\cdot {\rm sign}\, U=\psi(U)$. For anisotropic and
inhomogeneous media, we obtain the corresponding equation
(\ref{eqQUASILINEARW568CC}) below:

\medskip

{\bf Corollary 2.} {\it Let $D$ be a bounded domain in $\mathbb C$,
a scalar function $\sigma:D\to\mathbb R$ be in class
$L^{p^{\prime}}(D)$ for some $p^{\prime}>1$ and let $A:D\to\mathbb
S^{2\times 2}$ be a matrix function whose elements $a_{ij}(z)$,
$i,j=1,2$ are bounded and measurable.

Then there is a weak solution $u:D\to\mathbb R$ of class
$C^{\alpha}_{\rm loc}\cap W^{1,p}_{\rm loc}$ with $\alpha =1-2/p$
for some $p>2$ to the semi-linear Poisson type equation
\begin{equation} \label{SD1}
{\rm div}\, [A(z)\,\nabla\,u(z)]\ =\ \sigma(z)\cdot |u(z)|^{\lambda
-1} u(\xi)\ ,\ \ \ 0\ <\ \lambda\ <\ 1\  ,\ \ \ \ \ \ \ \mbox{a.e.
in $D$\ .}
\end{equation}}


Finally, we recall that in the combustion theory, see e.g.
\cite{Barenblat} and \cite{Po} and the references therein, the
following model equation
\begin{equation}\label{combustion}
{\partial u(z,t)\over \partial t}\ =\ {1\over \delta}\cdot \triangle
u\ +\ e^{u}\ ,\ \ \ \delta\ >\ 0\ ,\ t\geq 0,\ z\in D,
\end{equation}
takes a special part. Here $u\ge 0$ is the temperature of the
medium. We restrict ourselves here by the stationary case, although
our approach makes it possible to study the parabolic equation
(\ref{combustion}), see \cite{GNR}. The corresponding equation of
the type (\ref{SD}), see (\ref{eqQUASILINEARW568CCC}) below, appears
in anisotropic and inhomogeneous media with the function
$Q(u)=e^{-|u|}$ that is uniformly bounded at all.

\medskip

{\bf Corollary 3.} {\it Let $D$ be a bounded domain in $\mathbb C$,
a scalar function $\sigma:D\to\mathbb R$ be in class
$L^{p^{\prime}}(D)$ for some $p^{\prime}>1$ and let $A:D\to\mathbb
S^{2\times 2}$ be a matrix function whose elements $a_{ij}(z)$,
$i,j=1,2$ are bounded and measurable.

Then there is a weak solution $u:D\to\mathbb R$ of class
$C^{\alpha}_{\rm loc}\cap W^{1,p}_{\rm loc}$ with $\alpha =1-2/p$
for some $p>2$ to the semi-linear Poisson type equation
\begin{equation} \label{eqQUASILINEARW568CCC}
{\rm div}\, [A(z)\,\nabla\,u(z)]\ =\ \sigma(z)\cdot e^{-|u(z)|}\ \ \
\ \ \ \ \mbox{a.e. in $D$\ .}
\end{equation}}


{\bf Remark 6.} Such solutions $u$ in Corollaries 1--3 have the
representation as the composition $h\circ f^{\mu}|_D$, where
$f^{\mu}:\mathbb C\to\mathbb C$ is a $\mu -$conformal mapping in
Theorem B with $\mu$ extended onto $\mathbb C$ by zero outside of
$D$, and all $h$ are weak ge\-ne\-ra\-li\-zed harmonic functions
with sources of class $L_{p_*}(D_*)$, $D_*:=f^{\mu}(D)$ and
$p_*:=p^2/2(p-1)\in(2,p)$.

\section{On Hilbert problem for semi-linear Beltrami equations}

In this section, we prove results on existence, representation and
regularity of nonclassical solutions of the Hilbert boundary value
problem with arbitrary boun\-da\-ry data that are measurable with
respect to logarithmic capacity for semi-linear Beltrami equations
(\ref{SB}). On the history of the question and relevant definitions
as well notations can be found in Section 3 of our last preprint
\cite{GNRY*}.

\medskip

{\bf Theorem 3.}{\it\, Let $D$ be  a Jordan domain in $\mathbb C$
with the quasihyperbolic boundary condition, $\partial D$ have a
tangent q.e., $\lambda:\partial D\to\mathbb{C},\:
|\lambda(\zeta)|\equiv1$, be in $\mathcal{CBV}(\partial{D})$ and let
$\varphi:\partial D\to\mathbb{R}$ be measurable with respect to
logarithmic capacity.

Suppose also that $q:\mathbb C\to\mathbb C$ is a continuous function
with condition (\ref{q}), $\mu:D\to\mathbb C$ is of class
$L_{\infty}(D)$ with $k:=\|\mu\|_{\infty}<1$, $\mu$ is H\"older
continuous in an open neighborhood of $\partial D$ inside of $D$,
$\sigma: D\to\mathbb C$ has compact support in $D$, $\sigma\in
L_p(D)$ and condition (\ref{b}) holds for some $p>2$.

Then equation  (\ref{SB}) has a solution $\omega: D\to\mathbb C$ of
class $C^{\alpha}_{\rm loc}\cap W^{1,s}_{\rm loc}(D)$, where $\alpha
= 1-2/s$ and $s\in(2,p)$, that is smooth in the neighborhood of
$\partial D$ with the angular limits
\begin{equation}\label{eqLIMH} \lim\limits_{z\to\zeta , z\in D}\ \mathrm
{Re}\ \left\{\, \overline{\lambda(\zeta)}\cdot \omega(z)\, \right\}\
=\ \varphi(\zeta) \quad\quad\quad \mbox{q.e.\ on $\ \partial D$}\ .
\end{equation}}


{\bf Remark 7.} By the construction in the proof below, each such
solution has the representation $\omega =h\circ f|_D$, where
$f=f^{\mu}:\mathbb C\to\mathbb C$ is a $\mu -$conformal mapping from
Theorem B in \cite{GNRY*} with a suitable extension of $\mu$ onto
$\mathbb C$ and $h$ is a continuous solution of class
$W^{1,p_*}_{\rm loc}(D_*)$ for the semi-linear Vekua type equation
(\ref{SV}) with the multiplier $g\in L_{p_*}(D_*)$,
$p_*:=p^2/2(p-1)\in(2,p)$, as in (\ref{r}), $D_*:=f(D)$, which is a
generalized analytic function with a source in class $L_{p_*}(D_*)$
and has the angular limits
\begin{equation}\label{eqLIMH} \lim\limits_{w\to\xi , w\in D_*}\ \mathrm
{Re}\ \left\{\, \overline{\Lambda(\xi)}\cdot h(w)\, \right\}\ =\
\Phi(\xi) \quad\quad\quad \mbox{q.e.\ on $\ \partial D_*$}\ ,
\end{equation}
$\Lambda\,:=\,\lambda\circ f^{-1}|_{\partial D_*}$,
$\Phi\,:=\,\varphi\circ f^{-1}|_{\partial D_*}$. Also,
$s=p_*^2/2(p_*-1)\in(2,p_*)\subset(2,p)$.

\bigskip

{\bf Proof.} First of all, let us choose a suitable extension of
$\mu$ onto $\mathbb C$ outside of $D$. By hypotheses of Theorem 3
$\mu$ belongs to a class $C^{\alpha}$, $\alpha\in(0,1)$, for an open
neighborhood $U$ of $\partial D$ inside of $D$. By Lemma 1 in
\cite{GRY1} $\mu$ is extended to a H\"older continuous function
$\mu:U\cup\mathbb C\setminus D\to\mathbb{C}$ of the class
$C^{\alpha}$. Then, for every $k_*\in(k,1)$, there is an open
neighborhood $V$ of $\partial D$ in $\mathbb C$, where
$\|\mu\|_{\infty}\le k_*$ and $\mu$ in $C^{\alpha}(V)$. Let us
choose $k_*\in(k,1)$ so close to $k$ that $k_*C_p<1$ and set
$\mu\equiv 0$ outside of $D\cup V$.

By Theorem B in \cite{GNRY*}, there is a $\mu -$conformal mapping
$f=f^{\mu}:\mathbb C\to{\mathbb{C}}$ a.e. sa\-tis\-fying the
Beltrami equation (\ref{1}) with the given extended complex
coefficient $\mu$ in $\mathbb C$. Note that the mapping $f$ has the
H\"older continuous first partial derivatives in $V$ with the same
order of the H\"older continuity as $\mu$, see e.g. \cite{Iw} and
also \cite{IwDis}. Moreover, its Jacobian
\begin{equation}\label{6.13} J(z)\ne 0\ \ \ \ \ \ \ \ \ \ \ \ \
\forall\ z\in V\ , \end{equation} see e.g. Theorem V.7.1 in
\cite{LV}. Hence $f^{-1}$ is also smooth in $V_*:=f(V)$, see e.g.
formulas I.C(3) in \cite{Alf}.

Now, the domain $D_*:=f(D)$ satisfies the boundary quasihyperbolic
condition because $D$ is so, see e.g. Lemma 3.20 in \cite{GM}.
Moreover, $\partial D_*$ has q.e. tangents, furthermore, the points
of $\partial D$ and $\partial D^*$ with tangents correspond each to
other in one-to-one manner because the mappings $f$ and $f^{-1}$ are
smooth there. It is evident that the function
$\Lambda\,:=\,\lambda\circ f^{-1}|_{\partial D_*}$ belongs to the
class $\mathcal{CBV}(\partial{D_*})$.

Let us also show that the function $\Phi\,:=\,\varphi\circ
f^{-1}|_{\partial D_*}$ is measurable with respect to logarithmic
capacity. Indeed, for each open set $\Omega\subseteq\mathbb C$,
$\Phi^{-1}(\Omega)=f\circ \varphi^{-1}(\Omega)$, where the set
$\varphi^{-1}(\Omega)$ is  measurable with respect to logarithmic
capacity. Thus, it suffices to see that $f(S)$ is measurable with
respect to logarithmic capacity whenever $S$ is  measurable with
respect to logarithmic capacity.

Note for this goal that the quasiconformal mapping $f$ is H\"older
continuous on the compact set $\partial D$ and, thus, $C(f(S))=0$
whenever $C(S)=0$, see e.g. Remark 2.1 in \cite{GM}. Moreover, it is
known that Borel sets and, in particular, compact and open sets are
measurable with respect to logarithmic capacity, see e.g. Remark 2.2
in \cite{GM}. In addition, by definition a $C-$measurable set is a
union of a sigma-compactum (union of a countable collection of
compact sets) and a set $S$ with $C(S)=0$, see again Section 2,
especially formulas (2.5), in \cite{GM}. Thus, to conclude that $f$
translates $C-$measurable sets into $C-$measurable sets, it remains
to note that the ho\-meo\-mor\-phism $f$ translates compact sets
into compact sets.

Next, by Remark 2 in \cite{GNRY*} the function $g:D_*\to\mathbb C$
in (\ref{r}) belongs to class $L_{p_*}(D_*)$, where
$p_*=p^2/2(p-1)\in(2,p)$. Thus, by Theorem 2 in \cite{GNRY+} there
is a continuous solution $h$ of equation (\ref{SV}) that is a
generalized analytic function with a source of class $L_{p_*}(D_*)$
and that has the angular limits (\ref{eqLIMH}). Note that $h\in
W^{1,p_*}_{\rm loc}(D_*)$, see e.g. Theorems 1.16 and 1.37 in
\cite{Ve}. Finally, by Lemma 1 the function $\omega :=h\circ
f^{\mu}$ is a solution of (\ref{SB}) in class $C^{\alpha}_{\rm
loc}\cap W^{1,s}_{\rm loc}(D)$, where $\alpha = 1-2/s$ and
$s:=p_*^2/2(p_*-1)\in(2,p_*)\subset(2,p)$.\hfill $\Box$


In particular case $\lambda\equiv 1$, we obtain the corresponding
consequence of Theorem 3 on the Dirichlet problem for the
semi--linear Beltrami equations (\ref{SB}).

\section{On Poincare problem for semi-linear equations}

In this section we study the solvability of the Poincare
boundary-value problem for semi-linear Poisson type equations of the
form (\ref{SD}) in anisotropic and inhomogeneous media.

Given a matrix function $A: D\to\mathbb S^{2\times 2}$ and a $\mu
-$conformal mapping $f^{\mu}:D\to\mathbb C$, $\|\mu\|_{\infty}<1$,
we have already seen in Lemma 1 of \cite{GNR-}, by direct
computation, that if a function $T$ and the entries of $A$ are
sufficiently smooth, then
\begin{equation}\label{6.10k}
\mbox{div}\,[A(z)\nabla\,(T(f^{\mu}(z)))]\ =\ J(z)
\triangle\,T(f^{\mu}(z))\ . \end{equation} In the case $T\in
W^{1,2}_{\rm loc}$, we understand the identity (\ref{6.10k}) in the
distributional sense, see Proposition 3.1 in \cite{GNR}, i.e., for
all $\psi\in C^1_0(D)$,
\begin{equation}\label{6.10kd}
\int\limits_D\langle A\nabla(T\circ f^{\mu}),\nabla\psi\rangle\
dm_z\ =\ \int\limits_DJ(z)\langle M^{-1}((\nabla T)\circ
f^{\mu}),\nabla\psi\rangle\ dm_z\ ,
\end{equation}
where $M$ is the Jacobian matrix of the mapping $f^{\mu}$ and $J$ is
its Jacobian.

\medskip

{\bf Theorem 4.}{\it\, Let $D$ be  a Jordan domain with the
quasihyperbolic boundary condition, $\partial D$ have a tangent
q.e., $\nu:\partial D\to\mathbb{C},\: |\nu(\zeta)|\equiv 1$, be of
$\mathcal{CBV}(\partial{D})$ and $\varphi:\partial D\to\mathbb{R}$
be measurable with respect to logarithmic capacity.

Suppose also that $A\in M^{2\times 2}(D)$ has entries in a class
$C^{\alpha}(D)$, $\alpha\in(0,1),$ $\Sigma: D\to\mathbb R$ is a
function of class $L_p(D)$, $p>2 $, with a compact support in $D$
and $Q:\mathbb R\to\mathbb R$ is a continuous function with
condition (\ref{Q}).

Then there is a weak solution $u:D\to\mathbb R$ of class
$C^{1,\gamma}_{\rm loc}\cap W^{2,p}_{\rm loc}$, $\gamma
=\min(\alpha,\beta)$, $\beta =1-2/p,$ of the equation (\ref{SD})
that has the angular limits of its derivatives in the directions
$\nu =\nu(\zeta)$, $\zeta\in\partial D$,
\begin{equation}\label{eqLIMIT} \lim\limits_{z\to\zeta , z\in D}\ \frac{\partial u}{\partial \nu}\ (z)\ =\
\varphi(\zeta) \quad\quad\quad \mbox{q.e.\ on $\ \partial D$}\
.\end{equation}}


Here $u$ is called a {\bf weak solution} of equation (\ref{SD}) if
\begin{equation}\label{W}
\int\limits_{D}\{\langle A(z)\nabla u(z),\nabla \psi\rangle\ +\
\Sigma(z)\,Q(u(z))\,\psi(z)\}\ dm_z\ =\ 0\ \ \ \ \forall\ \psi\in
C^1_0(D)\ .
\end{equation}


{\bf Remark 8.} By the construction in the proof below, such a
solution $u$ has the representation $u =U\circ f|_D$, where
$f=f^{\mu}:\mathbb C\to\mathbb C$ is a $\mu -$conformal mapping from
Theorem B in \cite{GNRY*} with a suitable extension of $\mu$ in
(\ref{6.10}) to $\mathbb C$, and $U$ is a solution of class
$C^{1,\beta}\cap W^{2,p}_{\rm loc}$, $\beta =(p-2)/p,$ for the
semi-linear Poisson equation (\ref{SP}) with the multiplier (where
$J$ is the Jacobian of $f$) :
\begin{equation}\label{SOURCE}
G\ :=\ \frac{\Sigma}{J}\circ f^{-1}\ ,\ \ \ G\in L_{p}(D_*)\ ,\
D_*:=f(D)\ ,\end{equation}  that is a generalized harmonic function
with a source of the same class $L_{p}(D_*)$, which has the angular
limits
\begin{equation}\label{eqLIMHharm} \lim\limits_{w\to\xi , w\in D_*}\ \frac{\partial U}{\partial {\cal N}}\ (w)\ =\
\Phi(\xi) \quad\quad\quad \mbox{q.e.\ on $\ \partial D_*$}\ ,
\end{equation} where
\begin{equation}\label{N}
{\cal N}(\xi)\ :=\ \left\{\frac{\partial f}{\partial
\nu}\cdot\left|\frac{\partial f}{\partial
\nu}\right|^{-1}\right\}\circ f^{-1}(\xi)\ ,\ \ \ \xi\in\partial
D_*\ ,
\end{equation}
and
\begin{equation}\label{F}
\Phi(\xi)\ :=\ \left\{\varphi\cdot\left|\frac{\partial f}{\partial
\nu}\right|^{-1}\right\}\circ f^{-1}(\xi)\ ,\ \ \ \xi\in\partial
D_*\ .
\end{equation}


{\bf Proof.} By the hypotheses of the theorem $\mu$ given by
(\ref{6.10}) belongs to a class $C^{\alpha}(D)$, $\alpha\in(0,1)$,
and by Lemma 1 in \cite{GRY1} $\mu$ is extended to a H\"older
continuous function $\mu:\mathbb C\to\mathbb{C}$ of the class
$C^{\alpha}$. Then, for every $k_*\in(k,1)$, there is an open
neighborhood $V$ of $\overline D$, where $\|\mu\|_{\infty}\le k_*$
and $\mu$ is of class $C^{\alpha}(V)$. We may assume that $V$ is
bounded and set $\mu\equiv 0$ in $\mathbb C\setminus V$.

Let $f=f^{\mu}:\mathbb C\to{\mathbb{C}}$ be the $\mu -$conformal
mapping from Theorem B in \cite{GNRY+} with the given extended
complex coefficient $\mu$ in $\mathbb C$. Note that the mapping $f$
has the H\"older continuous first partial derivatives in $V$ with
the same order of the H\"older continuity as $\mu$, see e.g.
\cite{Iw} and also \cite{IwDis}. Moreover, its Jacobian
\begin{equation}\label{6.13} J(z)\ =\ |f_z|^2\, -\, |f_{\bar z}|^2\ >\ 0\ \ \ \ \ \ \
\forall\ z\in V\ , \end{equation} see e.g. Theorem V.7.1 in
\cite{LV}. Hence $f^{-1}$ is also smooth in $V_*:=f(V)$, see e.g.
formulas I.C(3) in \cite{Alf}.

Now, the domain $D_*:=f(D)$ satisfies the boundary quasihyperbolic
condition because $D$ is so, see e.g. Lemma 3.20 in \cite{GM}.
Moreover, $\partial D_*$ has q.e. tangents, furthermore, the points
of $\partial D$ and $\partial D^*$ with tangents correspond each to
other in a one-to-one manner because the mappings $f$ and $f^{-1}$
are smooth. In addition, the function $\cal N$ in (\ref{N}) belongs
to the class $\mathcal{CBV}(\partial{D_*})$ because
$$
\frac{\partial f}{\partial\nu}\ =\ f_z\cdot\nu\ +\ f_{\bar
z}\cdot\overline{\nu}\ ,\ \ \ \ \ \ \left|\frac{\partial
f}{\partial\nu}\right| \ \geq\ |f_z|\, -\, |f_{\bar z}|\ >\ 0\ ,
$$
and $\Phi$ in (\ref{F}) is measurable with respect to logarithmic
capacity by repeating arguments in the proof to Theorem 3.

Next, the source $G:D_*\to\mathbb R$ in (\ref{SOURCE}) belongs to
class $L_{p}(D_*)$ because by (\ref{6.13}) the function $J^{-1}\circ
f^{-1}$ is continuous and, consequently, bounded on the compact set
$\overline{D_*}$, see also point (vi) of Theorem 5 in \cite{ABe} on
the replacement of variables in integrals. Thus, by Theorem 4 in
\cite{GNRY+} there is a solution of class $C^{1,\beta}_{\rm
loc}(D_*)\cap W^{2,p}_{\rm loc}(D_*)$, $\beta =(p-2)/p,$ of the
semi-linear Poisson equation (\ref{SP}) with the multiplier
(\ref{SOURCE}) that is a generalized harmonic function with a source
of the same class $L_{p}(D_*)$ and which has the angular limits
(\ref{eqLIMHharm}) q.e. on $\partial D_*$.

Note that the function $u:=U\circ f$ belongs to class $W^{2,p}_{\rm
loc}(D)$ because $f$ is a quasi-isometry in $D$ of class $C^1$, see
e.g. 1.1.7 in \cite{Maz}. Finally, by Proposition 3.1 in \cite{GNR}
the function $u$ gives the desired solution of the equation
(\ref{6.8}) because by Lemma 10 and the point (i) of Theorem 5 in
\cite{ABe}
$$
\frac{\partial u}{\partial\nu}\ =\ u_z\cdot \nu\ +\ u_{\bar
z}\cdot\overline\nu\ =\ \nu\cdot(U_w\circ f\cdot f_{z}\ +\ U_{\bar
w}\circ f\cdot\overline{f_{\bar z}})\ +\
\overline{\nu}\cdot(U_w\circ f\cdot f_{\bar z}\ +\ U_{\bar w}\circ
f\cdot\overline{f_{z}})
$$
$$
=\ U_w\circ f\cdot (\nu\cdot f_{z}\ +\ \overline{\nu}\cdot f_{\bar
z})\
 +\ U_{\overline w}\circ f\cdot (\nu\cdot \overline{f_{\bar z}}\ +\ \overline{\nu}\cdot
 \overline{f_{z}})\ =\ U_w\circ f\cdot \frac{\partial f}{\partial
 \nu}\ +\ U_{\bar w}\circ f\cdot\overline{\frac{\partial f}{\partial\nu}}
$$
$$
=\ \left(\ {\cal N}\cdot U_w\ +\ \overline{\cal N}\cdot U_{\bar w}\
\right)\circ f\cdot \left| \frac{\partial f}{\partial \nu}\right|\
=\ \frac{\partial U}{\partial{\cal N}}\circ f\cdot \left|
\frac{\partial f}{\partial \nu}\right|\ ,
$$
where the direction $\cal N$ is given by (\ref{N}). This solution
$u$ belongs to the class $C^{1,\gamma}_{\rm loc}(D)$, $\gamma
=\min(\alpha,\beta)$, because by the above calculations with $\nu
=1$ and $i$
$$
u_x\ =\ U_w\circ f\cdot f_x\ +\ U_{\bar w}\circ f\cdot
\overline{f_x}\ ,\ \ \ u_y\ =\ U_w\circ f\cdot f_y\ -\ U_{\bar
w}\circ f\cdot \overline{f_y}\ ,\ \ \ z=x+iy\ .
$$


\medskip

{\bf Remark 9.} We are able to say more in the case of $\mathrm
{Re}\ n(\zeta)\overline{\nu(\zeta)}>0$, where $n(\zeta)$ is  the
inner normal to $\partial D$ at the point $\zeta$. Indeed, the
latter magnitude is a scalar product of $n=n(\zeta)$ and $\nu
=\nu(\zeta)$ interpreted as vectors in $\mathbb R^2$ and it has the
geometric sense of projection of the vector $\nu$ into $n$. In view
of (\ref{eqLIMIT}), since the limit $\varphi(\zeta)$ is finite,
there is a finite limit $u(\zeta)$ of $u(z)$ as $z\to\zeta$ in $D$
along the straight line passing through the point $\zeta$ and being
parallel to the vector $\nu$ because along this line
\begin{equation}\label{eqDIFFERENCE} u(z)\ =\ u(z_0)\ -\ \int\limits_{0}\limits^{1}\
\frac{\partial u}{\partial \nu}\ (z_0+\tau (z-z_0))\ d\tau\
.\end{equation} Thus, at each point with condition (\ref{eqLIMIT}),
there is the directional derivative
\begin{equation}\label{eqPOSITIVE}
\frac{\partial u}{\partial \nu}\ (\zeta)\ :=\ \lim_{t\to 0}\
\frac{u(\zeta+t\cdot\nu)-u(\zeta)}{t}\ =\ \varphi(\zeta)\ .
\end{equation}

\bigskip

In particular case of the Neumann problem, $\mathrm {Re}\
n(\zeta)\overline{\nu(\zeta)}\equiv 1>0$, where $n=n(\zeta)$ denotes
the unit interior normal to $\partial D$ at the point $\zeta$, and
we have by Theorem 4 and Remark 9 the following significant result.

\bigskip

{\bf Corollary 4.} {\it Let $D$ be a Jordan domain in $\Bbb C$ with
the quasihyperbolic boundary condition, the unit inner normal
$n(\zeta)$, $\zeta\in\partial D$, belong to the class
$\mathcal{CBV}(\partial D)$ and $\varphi:\partial D\to\mathbb{R}$ be
measurable with respect to logarithmic capacity.

Suppose also that $A\in M^{2\times 2}(D)$ has entries in a class
$C^{\alpha}(D)$, $\alpha\in(0,1)$, $\Sigma: D\to\mathbb R$ is a
function of class $L_p(D)$, $p>2 $, with a compact support in $D$
and $Q:\mathbb R\to\mathbb R$ is a continuous function with
condition (\ref{Q}).


Then one can find a weak solution $u:D\to\mathbb R$ of class
$C^{1,\gamma}_{\rm loc}\cap W^{2,p}_{\rm loc}$ with $\gamma
=\min(\alpha,\beta)$, $\beta =1-2/p,$ of equation (\ref{SD}) such
that q.e. on $\partial D$ there exist:


1) the finite limit along the normal $n(\zeta)$
$$
u(\zeta)\ :=\ \lim\limits_{z\to\zeta}\ u(z)\ ,$$

2) the normal derivative
$$
\frac{\partial u}{\partial n}\, (\zeta)\ :=\ \lim_{t\to 0}\
\frac{u(\zeta+t\cdot n(\zeta))-u(\zeta)}{t}\ =\ \varphi(\zeta)\ ,
$$

3) the angular limit
$$ \lim_{z\to\zeta}\ \frac{\partial u}{\partial n}\, (z)\ =\
\frac{\partial u}{\partial n}\, (\zeta)\ .$$}


{\bf Remark 10.} In addition, such a solution $u$ has the
representation $u =U\circ f|_D$, where $f=f^{\mu}:\mathbb
C\to\mathbb C$ is the $\mu -$conformal mapping from Theorem B in
\cite{GNRY*} with a suitable extension of $\mu$ in (\ref{6.10}) onto
$\mathbb C$ outside of $D$, described in the proof of Theorem 4, and
$U$ is a weak solution of the class $C^{1,\beta}_{\rm loc}\cap
W^{2,p}_{\rm loc}$, $\beta =1-2/p$, of the semi-linear Poisson
equation (\ref{SP}) with the multiplier $G$ in (\ref{SOURCE}) that
is a generalized harmonic function with a source of the same class
$L_{p}(D_*)$, which has the angular limits
\begin{equation}\label{neqLIMHharm}
\lim\limits_{w\to\xi , w\in D_*}\ \frac{\partial U}{\partial {n_*}}\
(w)\ =\ \varphi_*(\xi) \quad\quad\quad \mbox{q.e.\ on $\
\partial D_*$}\ ,
\end{equation} with
\begin{equation}\label{nN}
{n_*}(\xi)\ :=\ \left\{\frac{\partial f}{\partial
n}\cdot\left|\frac{\partial f}{\partial n}\right|^{-1}\right\}\circ
f^{-1}(\xi)\ ,\ \ \ \xi\in\partial D_*\ ,
\end{equation}
and
\begin{equation}\label{nF}
\varphi_*(\xi)\ :=\ \left\{\varphi\cdot\left|\frac{\partial
f}{\partial n}\right|^{-1}\right\}\circ f^{-1}(\xi)\ ,\ \ \
\xi\in\partial D_*\ .
\end{equation}

\section{The Poincare problem in physical applications}

By the discussion in Section 5, Theorem 4 on the Poincare
boundary-value problem for the Poisson type equations (\ref{SD})
gives the following consequences to mathematical models in
anisotropic and inhomogeneous media of diffusion with absorption,
plasma states and stationary burning, correspondingly.


{\bf Corollary 5.} {\it Let $D$ be  a Jordan domain with the
quasihyperbolic boundary condition, $\partial D$ have a tangent
q.e., $\nu:\partial D\to\mathbb{C},\: |\nu(\zeta)|\equiv 1$, be of
$\mathcal{CBV}(\partial{D})$ and $\varphi:\partial D\to\mathbb{R}$
be measurable with respect to logarithmic capacity.

Suppose also that $A\in M^{2\times 2}(D)$ has entries in a class
$C^{\alpha}(D)$, $\alpha\in(0,1)$, $\sigma: D\to\mathbb R$ is a
function of class $L_p(D)$, $p>2 $, with a compact support in $D$
and $Q:\mathbb R\to\mathbb R$ is a continuous function with
condition (\ref{Q}).

Then there is a weak solution $u:D\to\mathbb R$ of class
$C^{1,\gamma}_{\rm loc}\cap W^{2,p}_{\rm loc}$, $\gamma
=\min(\alpha,\beta)$, $\beta =1-2/p,$ of the semi-linear Poisson
type equation
\begin{equation} \label{eqQUASILINEARW568C}
{\rm div}\, [A(z)\,\nabla\, u(z)]\ =\ \sigma(z)\cdot u^{\lambda}(z)\
,\ \ \ 0\ <\ \lambda\ <\ 1\  ,\ \ \ \ \ \ \ \mbox{a.e. in $D$}
\end{equation}
satisfying the Poincare boundary condition on directional
derivatives
\begin{equation}\label{eqLIMIT68C} \lim\limits_{z\to\zeta}\ \frac{\partial u}{\partial \nu}\ (z)\ =\
\varphi(\zeta) \quad\quad\quad \mbox{q.e.\ on $\ \partial
D$}\end{equation} in the sense of the angular limits. }


{\bf Corollary 6.} {\it Let $D$ be  a Jordan domain with the
quasihyperbolic boundary condition, $\partial D$ have a tangent
q.e., $\nu:\partial D\to\mathbb{C},\: |\nu(\zeta)|\equiv 1$, be of
$\mathcal{CBV}(\partial{D})$ and $\varphi:\partial D\to\mathbb{R}$
be measurable with respect to logarithmic capacity.

Suppose also that $A\in M^{2\times 2}(D)$ has entries in a class
$C^{\alpha}(D)$, $\alpha\in(0,1)$, $\sigma: D\to\mathbb R$ is a
function of class $L_p(D)$, $p>2 $, with a compact support in $D$
and $Q:\mathbb R\to\mathbb R$ is a continuous function with
condition (\ref{Q}).

Then there is a weak solution $u:D\to\mathbb R$ of class
$C^{1,\gamma}_{\rm loc}\cap W^{2,p}_{\rm loc}$, $\gamma
=\min(\alpha,\beta)$, $\beta =1-2/p,$ of the semi-linear Poisson
type equation
\begin{equation} \label{eqQUASILINEARW568CC}
{\rm div}\, [A(z)\,\nabla\,u(z)]\ =\ \sigma(z)\cdot |u(z)|^{\lambda
-1} u(\xi)\ ,\ \ \ 0\ <\ \lambda\ <\ 1\  ,\ \ \ \ \ \ \ \mbox{a.e.
in $D$}
\end{equation}
satisfying the Poincare boundary condition on directional
derivatives (\ref{eqLIMIT68C}) q.e.}

\medskip

{\bf Corollary 7.} {\it Let $D$ be  a Jordan domain with the
quasihyperbolic boundary condition, $\partial D$ have a tangent
q.e., $\nu:\partial D\to\mathbb{C},\: |\nu(\zeta)|\equiv 1$, be of
$\mathcal{CBV}(\partial{D})$ and $\varphi:\partial D\to\mathbb{R}$
be measurable with respect to logarithmic capacity.

Suppose also that $A\in M^{2\times 2}(D)$ has entries in a class
$C^{\alpha}(D)$, $\alpha\in(0,1)$, $\sigma: D\to\mathbb R$ is a
function of class $L_p(D)$, $p>2 $, with a compact support in $D$
and $Q:\mathbb R\to\mathbb R$ is a continuous function with
condition (\ref{Q}).

Then there is a weak solution $u:D\to\mathbb R$ of class
$C^{1,\gamma}_{\rm loc}\cap W^{2,p}_{\rm loc}$, $\gamma
=\min(\alpha,\beta)$, $\beta =1-2/p,$ of the semi-linear Poisson
type equation
\begin{equation} \label{eqQUASILINEARW568CCC}
{\rm div}\, [A(z)\,\nabla\,u(z)]\ =\ \sigma(z)\cdot e^{-|u(z)|}\ \ \
\ \ \ \ \mbox{a.e. in $D$}
\end{equation}
satisfying the Poincare boundary condition on directional
derivatives (\ref{eqLIMIT68C}) q.e.}

\medskip

{\bf Remark 11.} Such solutions $u$ in Corollaries 5-7 have the
representation $u =U\circ f|_D$, where $f=f^{\mu}:\mathbb
C\to\mathbb C$ is a $\mu -$conformal mapping in Theorem B from
\cite{GNRY*} with a suitable extension of $\mu$ in (\ref{6.10}) to
$\mathbb C$ described in the proof of Theorem 4, and $U$ is a
solution of class $C^{1,\beta}\cap W^{2,p}_{\rm loc}$, $\beta
=(p-2)/p,$ for the semi-linear Poisson equation (\ref{SP}) with the
functions $Q(t)=t^{\lambda}$, $|t|^{\lambda -1}t$, $\lambda\in(0,1)$
and $e^{-|t|}$, correspondingly, and with the multiplier (here $J$
is the Jacobian of $f$) :
\begin{equation}\label{SOURCEs}
G\ :=\ \frac{\sigma}{J}\circ f^{-1}\ ,\ \ \ G\in L_{p}(D_*)\ ,\
D_*:=f(D)\ ,\end{equation} that is a generalized harmonic function
with a source $g$ of the same class $L_{p}(D_*)$, which satisfy the
Poincare boundary condition on directional derivatives
(\ref{eqLIMHharm}) in the sense of the angular limits q.e. on
$\partial D_*$.

\section{Neumann problem in physical applications}

In turn, Corollary 4 can be applied to the study of the physical
phenomena discussed by us in the last section. In this connection,
the particular cases of the function $Q(t)=t^{\lambda}$,
$|t|^{\lambda -1}t$, $\lambda\in(0,1)$, and $e^{-|t|}$ will be again
useful.


\newpage

{\bf Corollary 8.} {\it Let $D$ be a Jordan domain in $\Bbb C$ with
the quasihyperbolic boundary condition, the unit inner normal
$n(\zeta)$, $\zeta\in\partial D$, belong to the class
$\mathcal{CBV}(\partial D)$ and $\varphi:\partial D\to\mathbb{R}$ be
measurable with respect to logarithmic capacity.

Suppose also that $A\in M^{2\times 2}(D)$ has entries in a class
$C^{\alpha}(D)$, $\alpha\in(0,1)$, $\sigma: D\to\mathbb R$ is a
function of class $L_p(D)$, $p>2 $, with a compact support in $D$.

\medskip

Then one can find a weak solution $u:D\to\mathbb R$ of class
$C^{1,\gamma}_{\rm loc}\cap W^{2,p}_{\rm loc}$ with $\gamma
=\min(\alpha,\beta)$ and $\beta =1-2/p$  of the semi-linear Poisson
type equation (\ref{eqQUASILINEARW568C}) such that q.e. on $\partial
D$ there exist:

\bigskip

1) the finite limit along the normal $n(\zeta)$
$$
u(\zeta)\ :=\ \lim\limits_{z\to\zeta}\ u(z)\ ,$$

2) the normal derivative
$$
\frac{\partial u}{\partial n}\, (\zeta)\ :=\ \lim_{t\to 0}\
\frac{u(\zeta+t\cdot n(\zeta))-u(\zeta)}{t}\ =\ \varphi(\zeta)\ ,
$$

3) the angular limit
$$ \lim_{z\to\zeta}\ \frac{\partial u}{\partial n}\, (z)\ =\
\frac{\partial u}{\partial n}\, (\zeta)\ .$$}


{\bf Corollary 9.} {\it Under hypotheses of Corollary 8, there is a
weak solution $u$ $u:D\to\mathbb R$ of class $C^{1,\gamma}_{\rm
loc}\cap W^{2,p}_{\rm loc}$ with $\gamma =\min(\alpha,\beta)$ and
$\beta =1-2/p$  of the semi-linear Poisson type equation (\ref{SD})
such that q.e. on $\partial D$ all the conclusion 1)-3) of Corollary
8 hold, i.e., $u$ is a generalized solution of the Neumann problem
for (\ref{eqQUASILINEARW568CC}) in the given sense.}

\medskip

{\bf Corollary 10.} {\it Under hypotheses of Corollary 8, there is a
weak solution $u$ $u:D\to\mathbb R$ of class $C^{1,\gamma}_{\rm
loc}\cap W^{2,p}_{\rm loc}$ with $\gamma =\min(\alpha,\beta)$ and
$\beta =1-2/p$  of the semi-linear Poisson type equation (\ref{SD})
such that q.e. on $\partial D$ all the conclusion 1)-3) of Corollary
8 hold, i.e., $u$ is a generalized solution of the Neumann problem
for (\ref{eqQUASILINEARW568CCC}) in the given sense.}

\medskip

{\bf Remark 12.} Such solutions $u$ in Corollaries 6-10 have the
representation $u =U\circ f|_D$, where $f=f^{\mu}:\mathbb
C\to\mathbb C$ is a $\mu -$conformal mapping in Theorem B from
\cite{GNRY*} with a suitable extension of $\mu$ in (\ref{6.10}) to
$\mathbb C$ described in the proof of Theorem 4, and $U$ is a
solution of class $C^{1,\beta}\cap W^{2,p}_{\rm loc}$, $\beta
=(p-2)/p,$ for the semi-linear Poisson equation (\ref{SP}) with the
functions $Q(t)=t^{\lambda}$, $|t|^{\lambda -1}t$, $\lambda\in(0,1)$
and $e^{-|t|}$, correspondingly, and the multiplier $G$ in
(\ref{SOURCEs}), that is a generalized harmonic function with a
source $g$ of the same class $L_{p}(D_*)$, which satisfy the Neumann
boundary condition (\ref{neqLIMHharm}) in the sense of the angular
limits q.e. on $\partial D_*$.


{\bf \noindent Vladimir Gutlyanskii, Olga Nesmelova, Vladimir Ryazanov} \\
Institute of Applied Mathematics and Mechanics\\ of National Academy
of Sciences of Ukraine, Slavyansk,\\
\noindent{Bogdan Khmelnytsky National University of Cherkasy,\\
Physics Dept., Lab. of Math. Phys., Cherkasy, UKRAINE}\\
vgutlyanskii@gmail.com, star-o@ukr.net, vl.ryazanov1@gmail.com

\bigskip

\noindent {\bf Eduard Yakubov}\\
Holon Institute of Technology, Holon, ISRAEL,\\
yakubov@hit.ac.il, eduardyakubov@gmail.com

\end{document}